\documentclass{amsart}

\usepackage{amssymb} \usepackage{amsfonts} \usepackage{amsmath}
\usepackage{amsthm} \usepackage{epsfig} 

\newcommand{\mettifig}[1]{\epsfig{file=#1}}
\newcommand{\Sol}{{\rm Sol}}

\newtheorem{lemma}{Lemma}[section]
\newtheorem{teo}[lemma]{Theorem}
\newtheorem{prop}[lemma]{Proposition}
\newtheorem{cor}[lemma]{Corollary}

\theoremstyle{definition}

\newtheorem{quest}[lemma]{Question}

\theoremstyle{remark}
\newtheorem{rem}[lemma]{Remark}

\newcommand{\matN}{\ensuremath {\mathbb{N}}}
\newcommand{\matR} {\ensuremath {\mathbb{R}}}
\newcommand{\matQ} {\ensuremath {\mathbb{Q}}}
\newcommand{\matZ} {\ensuremath {\mathbb{Z}}}

\newcommand{\matCP} {\ensuremath {\mathbb{CP}}}

\newcommand{\calS} {\ensuremath {\mathcal{S}}}
\newcommand{\calM} {\ensuremath {\mathcal{M}}}

\def\interior#1{{\rm int}(#1)}

\newcommand{\JSJ}{{\rm JSJ}}
\newcommand{\Vol}{{\rm Vol}}

\newcommand{\timtil}{\begin{picture}(12,12)
\put(2,0){$\times$}\put(2,4.5){$\sim$}\end{picture}}

\newcommand{\finedimo}{{\hfill\hbox{$\square$}\vspace{2pt}}}

\author{Bruno Martelli}

\address{Dipartimento di Matematica \\%
Universit\`a di Pisa \\%
Via F.~Buonarroti 2 \\%
56127 Pisa, Italy%
}

\email{martelli@dm.unipi.it}

\urladdr{http://www.dm.unipi.it/$\sim$martelli/}

\title{Links, two-handles, and four-manifolds}



\thanks{The author is partially supported by the INTAS 
project ``CalcoMet-GT'' 03-51-3663}

\begin{document}

\begin{abstract}
We show that only finitely many links in a closed $3$-manifold can have homeomorphic 
complements,
up to twists along discs and annuli. Using the same techniques,
we prove that by adding $2$-handles on the same link we get only finitely
many smooth cobordisms between two given closed $3$-manifolds.

As a consequence, there are only finitely many
closed smooth 4-manifolds constructed from some Kirby diagram
with bounded number of crossings, discs, and strands,
or from some Turaev special shadow with bounded number of vertices.
(These are the $4$-dimensional analogues of Heegaard diagrams and special spines
for $3$-manifolds.)
We therefore get two filtrations 
on the set of all closed smooth orientable $4$-manifolds with finite sets.
The two filtrations are equivalent after linear rescalings, and
their cardinality grows at least as $n^{c \cdot n}$.
\end{abstract}

\maketitle

\section*{Introduction}
The main result of this paper can be interpreted as an extension of Thurston's Dehn
filling Theorem for hyperbolic $3$-manifolds to any $3$-manifold bounded by tori. 
The statement needs some preliminary definitions, and is therefore
deferred to Section~\ref{Dehn:surgery:section}.

Two independent finiteness results are consequences
of that theorem. The first one,
Theorem~\ref{link:complements:teo}, concerns the following
question: how many links can have homeomorphic complements in a closed $3$-manifold?
Concerning knots, a famous result
of Gordon and Luecke \cite{GL} states that the answer is ``one'' if we restrict
to knots in $S^3$. Bleiler, Hodgson, and Weeks have found two hyperbolic knots in a lens space 
with homeomorphic
complements \cite{BlHoWe}, but there are still no examples of distinct hyperbolic knots (in any
$3$-manifold) with the same \emph{oriented} complements. 

Gordon then extended his result from knots to a larger class of links in $S^3$ 
by showing in \cite{G} 
that at most $k!(38)^{k-1}$ such links with $k$ components can have homeomorphic
complements. If one considers all links in $S^3$, there are infinitely many 
having homeomorphic complements.
Theorem~\ref{link:complements:teo} below
states that those links are actually only finitely many up to some moves. The moves
consist of twisting along some discs and annuli.
This holds in any ambient closed $3$-manifold.

The second result, Theorem~\ref{main:teo}, 
is about $4$-manifolds obtained by adding $2$-handles on a given
link. As an application, it turns out that
only finitely many smooth closed $4$-manifolds can be described
from some Kirby diagram
with bounded numbers of crossings, discs, and strands,
or from some Turaev special shadow with a bounded number of 
vertices.

Kirby diagrams are the most common objects used to encode $4$-manifolds (and their
handle decompositions), see
\cite{GS} for a detailed overview. They are particularly efficient when the $4$-manifold
decomposes without $1$-handles, since they reduce to links in $S^3$ coloured with integers:
in that case, the manifold must be simply connected.
Turaev shadows can be simpler to handle for non-simply connected manifolds, or
for manifolds that need many $1$-handles in general. They are $2$-dimensional
polyhedra with half-integers on faces. They have proved recently very
effective in constructing $4$-manifolds with prescribed boundary: see the
work of Costantino and Thurston \cite{CoTh}.

Theorem~\ref{main:teo} implies that, by bounding the number of crossings, discs, and strands
in a Kirby diagram, or the number of vertices in a Turaev shadow, we get 
a filtration of the set of all smooth closed orientable $4$-manifolds (seen up to
diffeomorphism) into finite
sets. Theorem \ref{rescalings:teo} shows then that the two filtrations
are equivalent after linear rescalings. 

These finite filtrations allow one to study the enormous set of smooth 
$4$-manifolds from a probabilistic point of view, as Dunfield and Thurston recently
did with the set of $3$-manifolds \cite{DuTh}.
For instance, it makes sense in this context 
to calculate the probability for a smooth $4$-manifold to be
simply connected, complex, or symplectic (\emph{i.e.}~to admit one such structure). 
Theorem \ref{groups:teo}
shows that the total number of manifolds in the $n$-th set of the filtration
grows at least as $n^{cn}$. On the other hand, an easy bound on homology
together with well-known results of Freedman and Donaldson show
that the number of simply connected manifolds seen up to \emph{homeomorphism}
grows as $cn^2$. Many interesting questions then arise: 
does the number of smooth simply connected manifolds (up to diffeomorphism)
grow more than quadratic? Is the probability that a $4$-manifold has a definite
fundamental group zero?

Studying probabilities result in sending $n$ to infinity. 
On the other hand, it would be interesting to list the manifolds occurring for
low values of $n$, that is the ones lying in a first
segment of one filtration. One should compare the results with
analogous low-dimensional lists:
in the $3$-dimensional world, the most natural
filtration is probably given by Matveev complexity, which equals
(for irreducible manifolds) the minimum number of tetrahedra in a triangulation,
see for instance \cite{Mat}.
Manifolds with complexity $n\leqslant 10$ have been listed at various stages by 
Martelli, Matveev, and Petronio, 
and the results are summarized in \cite{survey}: 
all $\sim 700$ closed manifolds with complexity $\leqslant 8$ are graph
manifolds, while the first hyperbolic manifolds, appearing at complexity $9$, are the
$4$ smallest ones known. The number of hyperbolic manifolds then rapidly increases with $n$,
but the probability for a $3$-manifold to be hyperbolic is still unknown.

The first segment of one filtration introduced here 
has been recently studied by Costantino \cite{Co}:
the first $4$-manifolds are connected sums of simple manifolds as $\matCP^2$ and $S^2\times S^2$,
while more complicated manifolds seem to show up only for big values of $n$.
As in the $3$-dimensional case, it would be interesting to see what are the first
``more complicated manifolds'' in the list. For instance, what are the first simply connected
manifolds that do not decompose into $\matCP^2$ and $S^2\times S^2$ factors? Will
the $K3$ manifold be the first one? 

Of course, the same problems can be considered with respect to other finite filtrations,
like the one (which is probably the most natural)
given by the minimum number of $4$-simplexes in a triangulation, and can give
very different answers.
For the time being, it seems that Turaev shadows and Kirby diagrams can be studied with less
difficulties than triangulations (both from probabilistic and computational points of view), 
because of their tight relations with Dehn surgeries of $3$-manifolds and $2$-dimensional special
polyhedra: objects that have already been studied intensively by various authors.

The paper is organized as follows.
The results mentioned above are stated in detail in Section~\ref{finiteness:section},
except the core result of this paper, Theorem~\ref{convergence:teo}, 
which is stated and proved
in Section~\ref{Dehn:surgery:section}.
The results concerning links in $3$-manifolds are then proved 
in Section~\ref{links:section},
and those concerning $4$-manifolds in Section~\ref{handle:slides:section}.

\subsection*{Acknowledgements} 
The author would like to thank Fran\c cois Costantino for many helpful conversations
on shadows and $4$-manifolds. Trying with him to extend the natural properties of 
Matveev complexity of $3$-manifolds to the smooth $4$-dimensional category
has been the main motivation for this work.

\section{Two finiteness results} \label{finiteness:section}
We illustrate here the main consequences of Theorem~\ref{convergence:teo}.

\subsection{Link complements in $3$-manifolds}
Unlike knots~\cite{GL} in $S^3$, links in an arbitrary manifold $M$ 
are not determined by their complement. In fact, some moves may transform a given
link $L$ into another one, while preserving the complement: 
suppose $K$ is a component of $L$
with regular neighborhood $N(K)$, such that $M\setminus\interior{N(K)}$ contains an essential
disc $D$, transverse to the other components of $L$. Then one can make a full \emph{twist} 
of these components along $D$. Similarly, if the complement $M\setminus\interior{N(K_1\cup K_2)}$
of two components $K_1$ and $K_2$ of $L$ contains an essential annulus $A$ connecting the
two boundary tori, transverse to the other components of $L$, we can make a Dehn twist
along $A$. 

We prove here the following result.

\begin{teo}\label{link:complements:teo}
Let $M$ be a closed $3$-manifold. There are only finitely many links in $M$ sharing the same
complement, up to homeomorphisms of $M$ and twists along discs and annuli. 
\end{teo}

\begin{figure}
\begin{center}
\mettifig{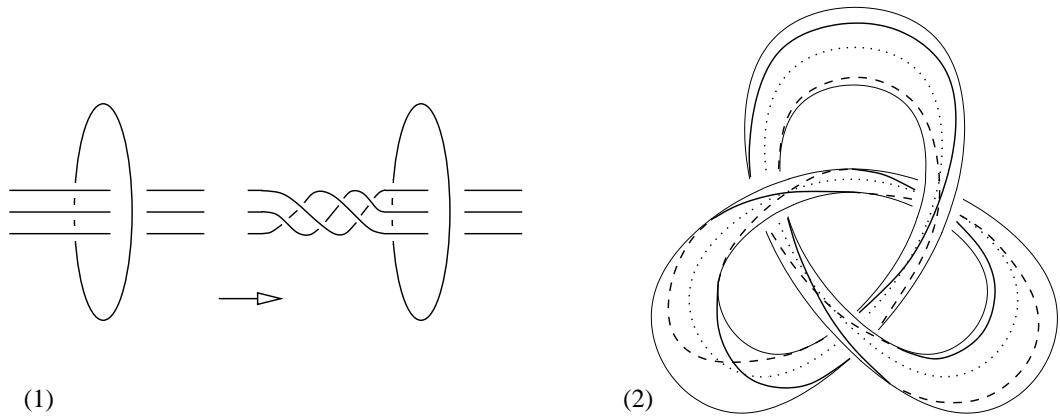, width = 14cm} 
\caption{A twist along a spanning disc (1) and two coaxial knots in $S^3$ (2): 
the core one is dotted.}
\label{twists:fig}
\end{center}
\end{figure}

The case $M=S^3$ has been studied by Gordon~\cite{G}. Inside $S^3$, one needs only to
consider discs that span unknotted components, 
as in Fig.~\ref{twists:fig}-(1), and annuli spanning two \emph{coaxial} components.
Two components $K_1$ and $K_2$ are coaxial when 
$K_1$ is an essential curve in the torus $\partial N(K_2)$, and is not a meridian of $N(K_2)$,
as in Fig.~\ref{twists:fig}-(2). 
The spanning 
annulus has one boundary component on $K_1$, and the other that may wind many times along $K_2$.

\begin{cor}\label{link:cor}
There are only finitely many links in $S^3$ sharing the same complement, up to 
twists along discs and annuli spanning unknotted and coaxial components.
\end{cor}

This generalizes Gordon's result~\cite{G}, saying that
there are only finitely many links in $S^3$ sharing the same complement, among those
not containing unknots and coaxial pairs. Actually, he proved that 
there are at most $k!(38)^{k-1}$
such links with $k$ components: it would be interesting to have an explicit 
bound also here.
Examples of $2$-component links with homeomorphic complements and not related by
twists along discs or annuli have been constructed by Berge~\cite{Be}.

\subsection{Handles and $4$-manifolds} \label{handles:subsection}
The set of all smooth 4-manifolds is a very difficult set to study.
A topological 
classification is impossible, essentially because every finitely presented group
is the fundamental group of some 4-manifold. 
Moreover, 
infinitely many non-diffeomorphic smooth manifolds can share the same topological
type~\cite{FM}.

Every smooth 4-manifold decomposes into handles. It turns out that
$0$-, $1$-, $3$-, and $4$-handles can be attached essentially
in a unique way~\cite{LP}, whereas the huge
variety of smooth 4-manifolds is due to the
many possibilities one has to attach 2-handles.

Two-handles are encoded via a framed link in the boundary
$3$-manifold. The following result
says roughly that the huge variety
of smooth 4-manifolds with fixed boundary is due to the variety of links 
in $3$-manifolds, not to their framings.

\begin{teo} \label{main:teo}
Let $M$ and $N$ be two closed $3$-manifolds, and $L\subset M$ be a link.
There are only finitely many smooth cobordisms (up to diffeomorphism)
between $M$ and $N$ obtained by adding $2$-handles to $M$ along $L$ (with some framing).
\end{teo}

\begin{figure}[t]
\begin{minipage}{.02\textwidth}\hfil
\end{minipage}
\begin{minipage}{.2\textwidth}
\centering
\mettifig{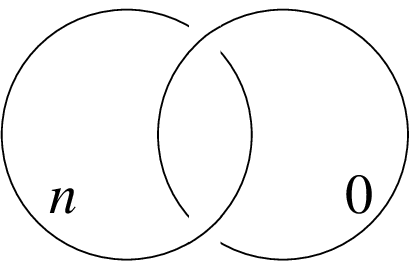, width = .9\textwidth}
\vglue -4mm
\end{minipage}
\begin{minipage}{.02\textwidth}\hfil
\end{minipage}
\begin{minipage}{.66\textwidth}
{Here $M = N=S^3$ for all $n\in\matZ$. 
The resulting cobordisms, capped with $2$ balls, give $S^2\times S^2$ for even $n$,
and $S^2\timtil S^2 \cong 
\matCP^2\#\overline{\matCP^2}$ for odd $n$~\cite{GS}.
}
\end{minipage}
\caption{Handles on the Hopf link. Framings are described via integers, 
as usual~\cite{K, GS}. The infinitely many cobordisms, 
parametrized by $n\in\matZ$, are actually only two up to diffeomorphism,
in accordance with Theorem~\ref{main:teo}.}
\label{Hopf:fig}
\end{figure}

To be precise, $2$-handles are attached 
to the manifold $M\times [0,1]$ along its boundary component 
$M\times\{1\}$, which we identify with $M$.
All cobordisms obtained by adding $2$-handles on $L$, with varying $N$,
are homotopically equivalent. But in the smooth category the $3$-manifold $N$ has to be fixed,
because in general infinitely many distinct $N$'s are obtained from a fixed $L\subset M$.
The same $N$ can nevertheless be obtained via infinitely many 
different framings on the same link, as Fig.~\ref{Hopf:fig} shows.

\subsection{$2$-polyhedra in $4$-manifolds}
Theorem~\ref{main:teo} can be used to derive other finiteness results
on $4$-manifolds. 
Let a \emph{$k$-skeleton} $P^k\subset M^n$ in a smooth $n$-manifold 
be a $k$-dimensional polyhedron such that the complement of its regular neighborhood
consists of handles of dimension higher than $k$.

If $k=n-1$ and $P^k$ has a nice local structure 
(\emph{i.e.} it is \emph{special}), the polyhedron $P^k$ determines $M^n$~\cite{Mat:alta:dim}.
A similar result in codimension-$2$ seems hopeless at a first sight: 
infinitely many
$3$-manifolds share the same $1$-skeleton, no matter how nice it is.
Quite surprisingly, things get better in dimension $4$: the nice local
structure needed here is local flatness.

A 2-dimensional polyhedron $P\subset M$ in a smooth $4$-manifold $M$ is 
\emph{locally flat}\footnote{The smooth and PL notions of local flatness coincide for
surfaces but differ heavily for polyhedra:
see Subsection~\ref{shadows:subsection} for a discussion on this point.
We use here the stronger smooth notion.}
if it is locally contained in smooth 3-discs in $M$.

\begin{teo}\label{polyhedra:teo}
A $2$-dimensional compact polyhedron is the locally flat $2$-skeleton
of only finitely many closed orientable smooth $4$-manifolds.
\end{teo}

\subsection{Finite filtrations on the set of smooth closed $4$-manifolds} \label{closed:subsection}
\begin{figure}
\begin{center}
\mettifig{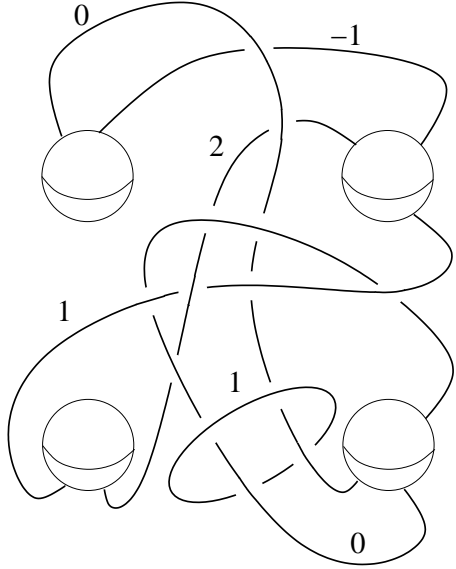, width = 6 cm} 
\caption{A Kirby diagram. It consists of some couples of discs and some strands, coloured
with integers (determining the framing). 
The strands are attached to each couple of discs symmetrically as shown here.}
\label{Kirby:fig}
\end{center}
\end{figure}

A handle decomposition of a closed orientable $4$-manifold 
can be encoded via a planar \emph{Kirby diagram}~\cite{K, GS},
which consists of some couples of discs (respresenting $2$-spheres in $S^3$ 
that encode $1$-handles) and a link diagram (encoding $2$-handles) 
with some strands having endpoints on
the discs as in Fig.~\ref{Kirby:fig}.
The strands are coloured with integers, determining their framings.
Higher handles need not to be encoded
thanks to the Laudenbach-Poenaru Theorem~\cite{LP}. 
We define the \emph{weight} of a Kirby diagram as the sum of the 
numbers of crossings, discs, and strands.
The following result is a corollary of Theorem~\ref{main:teo}.

\begin{cor}\label{Kirby:cor}
The number of closed orientable smooth $4$-manifolds described via a Kirby diagram
with weight at most $n$ is finite, for each $n$.
\end{cor}

\begin{rem}\label{filtration1:rem}
Let $U_n$ be the set of all closed orientable smooth manifolds described via some 
Kirby diagram
with weight at most $n$. We get a filtration 
$U_1\subset\ldots\subset U_n\subset\ldots$
on the set of all closed oriented smooth $4$-manifolds with finite sets $U_n$.
\end{rem}

\begin{figure}
\begin{center}
\mettifig{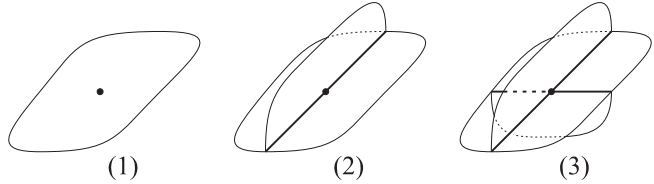} 
\caption{Neighborhoods of points in a special polyhedron.}
\label{new_standard_nhbds:fig}
\end{center}
\end{figure}

A compact 2-dimensional polyhedron $P$ is \emph{special} if every point
of $P$ has one of the regular neighborhoods shown in 
Fig.~\ref{new_standard_nhbds:fig}, and if 
the stratification given by the three
types of points gives a cellularization of $P$ (\emph{i.e.}~points
of type (1) form discs and points of type (2) form segments).
A \emph{vertex} is a point as in Fig.~\ref{new_standard_nhbds:fig}-(3).

Inspired by Turaev's terminology~\cite{Tu}, we
call \emph{special shadow} $P\subset M$ a locally flat special $2$-skeleton
in a closed orientable smooth $4$-manifold $M$.
Since there are only finitely
many special polyhedra with at most $n$ vertices~\cite{Mat}, 
we deduce the following from Theorem~\ref{polyhedra:teo}.

\begin{cor} \label{shadow:cor}
The number of closed orientable smooth $4$-manifolds having a special shadow
with at most $n$ vertices is finite, for each $n$.
\end{cor}

\begin{rem}\label{filtration2:rem}
Every closed orientable smooth $4$-manifold has a special shadow~\cite{Tu}.
Therefore, as in Remark~\ref{filtration1:rem}, by defining
$V_n$ as the set of 
manifolds having a shadow with $n$ vertices
at most, we get a filtration on all closed 
oriented smooth $4$-manifolds with finite sets.
\end{rem}

We make some comments on the two filtrations
introduced above in Remarks~\ref{filtration1:rem} and~\ref{filtration2:rem}.
The filtrations are equivalent after
linear rescalings, as the following result shows.

\begin{teo} \label{rescalings:teo}
We have $U_n\subset V_{3n}$ and $V_n\subset U_{9n+8}$ for all $n>0$.
\end{teo}

We now give some information on the cardinality of $U_n$ and $V_n$.
We say that a sequence $a_n$ of integers \emph{grows as $n^n$} if there are
constants $0<c<C$ such that $n^{c\cdot n}< a_n < n^{C\cdot n}$ for all $n\gg 0$.

\begin{teo} \label{groups:teo}
The number of distinct groups that are fundamental groups of manifolds
in $U_n$ grows as $n^n$. The same result holds for $V_n$.
\end{teo}
\begin{cor}
Each set $U_n$ or $V_n$ contains at least $n^{c\cdot n}$
distinct manifolds, for some $c>0$ and for all $n\gg 0$.
\end{cor}

The proof of Theorem~\ref{groups:teo} makes use of the techniques of~\cite{FriMaPe3},
that is of Mostow rigidity for hyperbolic $3$-manifolds with geodesic boundary.
Concerning simply connected manifolds, by constructing connected
sums of $\matCP^2$ one immediately gets:

\begin{prop} \label{simply:connected:prop}
There are at least $\frac 14n^2$ distinct simply connected manifolds in $U_n$, for all $n$.
\end{prop}

On the other hand, 
using the celebrated theorems of Freedman~\cite{F} and Donaldson~\cite{D} we easily show:
\begin{prop} \label{homeomorphic:prop}
There are at most $\frac 5{16} n^2$ simply connected 
manifolds in $U_n$ up to homeomorphism, for $n\gg 0$.
\end{prop}
\begin{quest} What is the growth type of the number of simply connected smooth manifolds
in $U_n$ (or $V_n$)? Is it bigger than quadratic?
\end{quest}

\section{Dehn filling} \label{Dehn:surgery:section}

In this Section, all $3$-manifolds
will be compact and orientable, and possibly with boundary consisting of tori.
We define here a volume and a Euler number for any 
geometrizable $3$-manifold, which
extend both the volume of a hyperbolic manifold and the Euler number of a Seifert manifold.
We then use these definitions to state and prove Theorem~\ref{convergence:teo}.

\subsection{Slopes}
Let $T$ be the torus. After fixing a homology basis $(m,l)$ for
$H_1(T;\matZ)$, every \emph{slope} on $T$ (\emph{i.e.}~isotopy class of simple closed
essential curves) is determined by its unsigned homology class $\pm(pm+ql)$,
thus by the number $q/p\in\matQ\cup\{\infty\}$
The \emph{distance} $\Delta(q/p, s/r)$ of two slopes
is their minimal geometric intersection, equal to $|ps-qr|$.
(This is not really a distance, since the triangle inequality does not hold.)
The set of all slopes in $T$ inherits the topology of
$\matQ\cup\{\infty\}$,
which does not depend on the chosen basis $(m,l)$.
Its completion $\matR\cup\{\infty\}\cong S^1$ can be seen as the space of all
geodesic foliations of $T$ (with some flat metric), up to isotopy.
We will often use below the following fact.

\begin{prop}\label{distance:prop}
Suppose $\{s_i\}_{i\in\matN}$ and $\{s_i'\}_{i\in\matN}$ are sequences of 
slopes such that 
\begin{itemize}
\item $s_i\to \lambda$ for some foliation $\lambda$ with $s_i\neq \lambda$ for all $i$;
\item the sequence $\{s_i'\}$ is contained in some closed set not containing $\lambda$. 
\end{itemize}
In this case we have $\Delta(s_i,s_i')\to\infty$.
\end{prop}
\begin{proof}
We can suppose no slope or foliation 
involved is $\infty$. Therefore $s_i=q_i/p_i$ with $|p_i|\to\infty$
and $|s_i - s_i'| > k$ for all $i$. We get:
$$\Delta(s_i,s_i')=\Delta(q_i/p_i,q'_i/p'_i)=|q_ip'_i-p_iq'_i|=|(s_i-s'_i)p_ip'_i|\to\infty.$$
\end{proof}

\subsection{Dehn filling}
Let $M$ be a manifold with $\partial M$ containing some tori $T_1,\ldots,T_k$.
We denote by $M(s)$ the \emph{Dehn filling} of $M$ along the vector of slopes
$s=(s_1,\ldots,s_k)$, obtained from $M$ by
attaching to $T_i$ a solid torus with a map sending the meridian 
to $s_i$, for all $i$. 

\subsection{Seifert manifolds}~\label{Seifert:subsection}
We recall the definition and some properties of Seifert manifolds.
Let $M$ be an oriented $S^1$-fibering over some compact surface $F$ with boundary. 
A homology basis $(m,l)$ is defined for each boundary torus $T\subset \partial M$
by taking $m$ as the boundary of a fixed section
and $l$ as the fiber.
If $F$ is orientable we orient it and orient coherently $m$ and $l$, otherwise
we choose any orientation. 
The Dehn filling $N = M(q_1/p_1,\ldots,q_k/p_k)$ fibers
over the orbifold $\Sigma$ obtained by capping $k$ boundary components of $F$ with
discs having points of cone angle $2\pi/p_1,\ldots,2\pi/p_k$ (we assume \emph{wlog} that
$p_i>0$ for all $i$).
Two important invariants of that fibration are
the Euler characteristic $\chi(\Sigma)$ of $\Sigma$ and the Euler number $e$, given by
$$\chi(\Sigma)=\chi(F)+\sum_{i=1}^k\frac 1{p_i}\qquad {\rm and} \qquad
e = \sum_{i=1}^k\frac {q_i}{p_i}.$$

When $M$ has boundary, $e$ is only defined up to additive integers, and we
require that $0\leqslant e <1$.
By substituting a pair $(q_i/p_i, q_j/p_j)$ with 
$(q_i/p_i+1, q_j/p_j-1)$ we get the same fibration.
The manifold $N$ is called \emph{Seifert}. 
The following result will be needed below.

\begin{prop}\label{small:Seifert:prop}
Let a closed $N$ fiber over $\Sigma$ with Euler number $e$: 
\begin{itemize}
\item
if $\chi(\Sigma)\leqslant 0$, the non-negative number $|e|$ depends only on $N$;
\item
if $\chi(\Sigma)>0$, we have $|\pi_1(N)|\geqslant |e|$.
\end{itemize}
\end{prop}
\begin{proof}
If $\chi(\Sigma) \leqslant 0$, then $N$ has a unique fibration, 
except some flat cases~\cite{Sco}
where $\chi = e = 0$, and we are done.
Suppose $\chi > 0$ and $|\pi_1(N)|<\infty$. Then the fibration of $N$ lifts to a
fibration of $S^3$ with Euler number $\tilde e$.
A fibration of $S^3$ has base space $S^2$ and at most two cone points.
Then it is constructed by filling $A\times S^1$ with some 
$(q_1/p_1, q_2/p_2)$, where $A$ is the annulus.
Since the total space is $S^3$, we have $|p_1q_2+p_2q_1|=1$. 
Therefore $|\tilde e| = 1/|p_1p_2|\leqslant 1$.
By \cite[Theorem 3.6]{Sco} we have 
$|e| \leqslant |\pi_1(N)|\cdot|\tilde e| \leqslant |\pi_1(N)|$.
\end{proof}

\subsection{Sol-manifolds} 
A $3$-manifold fibering over a $1$-orbifold with fibers consisting
of tori and Klein bottles, which does not admit a Seifert fibration, has a 
Sol-geometry~\cite{Sco}. It consists of a torus fibering over $S^1$
or of two twisted interval bundles over a Klein bottle glued along their boundaries.

\subsection{Geometric decomposition of a 3-manifold}
An irreducible orientable compact 3-manifold, possibly with boundary
consisting of tori, is \emph{geometrizable} if it
satisfies Thurston's Geometrization Conjecture.
A geometrizable 3-manifold has a
unique \emph{geometric decomposition} along
embedded tori and Klein bottles into blocks having one of the
$8$ three-dimensional geometries,
constructed from the set of tori of the \JSJ\ 
decomposition\footnote{Here, the \JSJ\ decomposition of a manifold of type \Sol\ 
is trivial.},
by substituting each torus bounding a 
twisted interval bundle over a Klein bottle $K$ with the core $K$.
A non-empty decomposition can be easily checked to be geometric, as the following shows.
\begin{prop}\label{geometric:prop}
Let an orientable $M$ be decomposed along a non-empty set of tori 
and Klein bottles into some Seifert or hyperbolic blocks.
Such a decomposition is the geometric one if and only if the following holds:
\begin{enumerate}
\item
every Seifert block fibers over an orbifold $\Sigma$ with $\chi(\Sigma)<0$;
\item
the fibrations of the blocks adjacent to a torus or a Klein bottle $S$ do not induce
the same fibration on $S$.
\end{enumerate}
\end{prop}

When $S$ is a Klein bottle, we mean that the fibration of the 
single block adjacent to $S$ should not
induce a fibration on $S$ (a Klein bottle admits only two non-isotopic
fibrations). A Seifert manifold with $\chi<0$ has a unique fibration~\cite{Sco}.

\subsection{Generalized volume}
Let $M$ be a geometrizable manifold.
We now define a non-negative quantity $\Vol(M)$ which generalizes the volume of a hyperbolic
manifold. If $M$ is of type \Sol\ or Seifert with $\chi\geqslant 0$, 
we set $\Vol(M)=0$.
Otherwise, we define $\Vol(M)$ as the sum
of the volumes of the hyperbolic blocks, minus the sum of the $\chi$'s of the Seifert blocks
in its geometric decomposition.

\subsection{Generalized Euler number}
Let $M$ be a geometrizable manifold.
We now define a non-negative quantity $e(M)$ which generalizes 
the Euler number of a fibration, as follows:
\begin{itemize}
\item if $M$ admits a Seifert fibration with $\chi\leqslant 0$ and some Euler number $e$, 
then $e(M)=|e|$ is well-defined by Proposition~\ref{small:Seifert:prop};
\item if $M$ admits a Seifert fibration with $\chi>0$ we set $e(M)=|\pi_1(M)|$ 
when $|\pi_1(M)|<\infty$ and $e(M)=0$ otherwise;
\item if $M$ is hyperbolic we set $e(M)=0$;
\item if $M$ is a torus fibering over $S^1$ with monodromy $\psi$, we set 
$$e(M)=\min\{\Delta(m,\psi(m))+\Delta(l,\psi(l))\}$$ among all pairs of slopes
$(m,l)$ in the fiber with $\Delta(m,l)=1$;
\item if $M$ is the union of two twisted interval bundles $N$ and $N'$
over the Klein bottle, take $e(M)=\max\{\Delta(l,l')\}$ among all fibrations
of $N$ and $N'$ with fibers $l\subset\partial N$ and $l'\subset\partial N'$.
\end{itemize}
Otherwise, $M$ has a non-trivial geometric decomposition into blocks, and we want
$e(M)$ to measure how complicated the gluing maps between them are. To do this,
we start by defining for each block $N$ 
and each abstract boundary torus $T\subset\partial N$ a finite set
of \emph{preferred slopes}, which contains at least two elements and depends
only on $N$ and $T$.

If $N$ is hyperbolic, we define a preferred slope
to be any slope having the shortest or second
shortest length in one cusp section.
If $N$ is Seifert, the situation is more complicated because the fiber is the
only slope which is intrinsically defined.
We see $N$ as the filling of some $S^1$-bundle over a surface $F$ along
some slopes $(q_1/p_1,\ldots,q_k/p_k)$ with $0<q_i/p_i<1$.
We define two preferred slopes on $T$: the fiber of $N$ and $\partial F\cap T$.
Although the first one is intrinsic, the second
one depends on the \emph{section} $F$.

Let $S$ be a surface of the geometric decomposition. If $S$ is a torus,
we define $\Delta_S$ as the maximum of $\Delta(s_1,s_2)$, where
$s_1$ and $s_2$ are some preferred slopes on the
two adjacent blocks. If $S$ is a Klein bottle,
a small neighborhood $W$ admits two fibrations, and
we define $\Delta_S$ as the maximum of $\Delta(s_1,s_2)$
on the torus $\partial W$,
where $s_1$ is some preferred slope of the adjacent block,
and $s_2$ is one of the two fibers of $W$.

We now define $\Delta$ as the maximum of $\Delta_S$ as $S$ varies.
This quantity depends on the sections $F$
chosen. Finally, we define $e(M)$ as the minimum of $\Delta$
as the sections vary.

\subsection{Properties of the invariants}
The following result is not
used elsewhere in this paper.
\begin{teo}
Let $\calM$ be the set of all geometrizable $3$-manifolds:
\begin{enumerate}
\item 
the set $\Vol(\calM)\subset\matR$ is well-ordered; 
\item 
there are only finitely many manifolds $M$ with $\Vol(M) = v_0$ and $e(M)\leqslant e_0$,
for all $v_0$ and $e_0$.
\end{enumerate}
\end{teo}
\begin{proof}
The set $\Vol_H$ of volumes of hyperbolic manifolds is well-ordered~\cite{bibbia}. 
The set $\Vol_S$ of volumes of Seifert manifolds is
$$\Vol_S = \{n - (1/p_1+\cdots+ 1/p_k) \ \big| \ n>0,\ 0\leqslant k\leqslant n+2\}\cap\matR_+,$$
which is well-ordered. Therefore
$$\Vol(\calM) = \{x_1+\cdots+x_k\ \big| \ x_i\in \Vol_H\cup \Vol_S\}$$
is also well-ordered. 

We now prove point (2). 
There are only finitely many hyperbolic manifolds
with volume $v_0$~\cite{BP}. 
There are also only finitely many $2$-orbifolds with given $\chi$.
On each orbifold $\Sigma$, there are only finitely many Seifert fibrations with Euler number 
smaller than $e_0$. Therefore point (2)
holds for hyperbolic and Seifert manifolds.

We note the following: let $T_1$ and $T_2$
be two tori, each $T_i$ equipped with a homology basis $(m_i, l_i)$. For each $K>0$ there are
only finitely many homeomorphisms $\psi:T_1\to T_2$ (up to isotopy) with 
$\Delta(\psi(m_1),m_2)\leqslant K$ and $\Delta(\psi(l_1),l_2)\leqslant K$.
Therefore there are only finitely many \Sol-manifolds with $e = e_0$, and
the case $v_0=0$ is done, since there are also only finitely many elliptic manifolds 
with bounded $|\pi_1|$.

We are therefore left with the case where $M$ is not geometric.
The manifold $M$ decomposes into some geometric pieces of volumes $v_1,\ldots, v_k$.
Since $\Vol(\calM)$ is well-ordered, there are only finitely many possible sequences $v_1,\ldots, v_k$
with $v_1+\cdots +v_k = \Vol(M)$. 
We have just proved that there are only finitely many geometric
manifolds with volume $v_i$. 
Therefore there are only finitely many possible pieces in the decomposition.
Finally, the condition $e(M)=e_0$ ensures that
only finitely many gluings are admitted on every torus, and we are done.
\end{proof}

\subsection{Sequences of slopes}
Let $N$ be a 3-manifold bounded by $k$ tori. The vectors of slopes on $\partial N$
form a set $(\matQ\cup\{\infty\})^k$, whose completion is a $k$-torus $T^k$, which
has a geometric interpretation as the set of geodesic foliations on $\partial N$.

A filled closed manifold $N(s)$ is associated to every rational point
$s=(s_1,\ldots,s_k)$ of $T^k$. Theorem~\ref{convergence:teo} below 
gives an answer to the following question:
if $\{s^i\}_{i\in\matN}$ is a sequence of rational points on $T^k$, what can we say
about the sequence of manifolds $\big\{N(s^i)\big\}$?

Since $T^k$ is compact, a subsequence of $\{s^i\}$ converges to a foliation $\lambda$.
Set $s^i = (s^i_1,\ldots,s^i_k)$ and $\lambda = (\lambda_1,\ldots,\lambda_k)$.
We call a sequence \emph{essential} if $s^i_j\neq \lambda_j$ for all $i$.
If $\{s^i\}$ contains no essential subsequence, it contains one
having $s^i_j = \lambda_j$ constantly for
all $i$, which can be studied on the manifold $N(\lambda_j)$ having
one boundary component less. For this reason we 
restrict our attention to essential sequences.

This discussion motivates the following result.
For a sequence $\{a_i\}_{i\in\matN}$ of real numbers, we write $a_i\nearrow a$ when
$a_i<a$ for all $i$ and $a_i\to a$. 

\begin{teo}\label{convergence:teo}
Let $N$ be an irreducible compact $3$-manifold bounded by $k$ tori.
Let $\{s^i\}_{i\in\matN}$ be an essential sequence of (vectors of) slopes, converging to 
some foliation $\lambda = (\lambda_1,\ldots,\lambda_k)$. 

Suppose no $\lambda_j$ is a slope that bounds a disc, 
and no distinct $\lambda_j$ and $\lambda_{j'}$ are slopes that cobound an annulus
in $N$.
After passing to a subsequence,
every $N(s^i)$ is geometrizable and one of the following holds:
$$\Vol(N(s^i))\nearrow \Vol(N) \qquad {\rm or} \qquad e(N(s^i))\nearrow\infty.$$
\end{teo}

\begin{rem}
Theorem~\ref{convergence:teo} is false without the hypothesis on discs and annuli.
For instance, if $N$ is a solid torus and $s^i$ is a slope intersecting the meridian
$\lambda$ once and winding $i$ times along it, we get $s^i\to \lambda$ 
and $N(s^i)=S^3$ not depending on $i$. Analogously, if $N$ is a Seifert manifold
with two boundary components and $s^i = (q/p + i, q'/p' - i)$, 
then the limit slopes $(\infty, \infty)$ bound a fibred annulus 
and $N(s^i)$ does not depend on $i$.
\end{rem}

\subsection{Examples}
We show some of the phenomena that can occur. 
Recall that, with our choices of homology bases in the boundary 
of a Seifert manifold, the coordinate of the fiber is $\infty$, and a $(q/p)$-Dehn filling
produces a new singular fiber if and only if it is not an integer.
\begin{itemize}
\item
Let $N$ be the solid torus and the limit $\lambda$ be any foliation distinct
from the meridian: we get $e(N(s^i)) = |\pi_1(N(s^i))|\nearrow\infty$;
\item
suppose $N$ is not a solid torus or an interval bundle, and $\lambda_j\neq\infty$ 
on every Seifert block of $N$:
then (for big $i\in\matN$) the geometric decomposition of $N(s^i)$ is
induced by the one of $N$ and $\Vol(N(s^i))\nearrow\Vol(N)$;
\item
let $N$ be a Seifert manifold with $\chi<0$
and one boundary component, filled with an integer $s^i=i \nearrow \infty$:
we get $\Vol(M(s^i))=\Vol(M)$ and $e(M(s^i))\nearrow\infty$;
\item
let $P$ be the pair-of-pants, and let $N$ decompose into $P\times S^1$ and
two hyperbolic blocks $H$ and $H'$ with one cusp; then $\partial N$ is one torus, and:
\begin{itemize}
\item if $s^i=i$, each $N(s^i)$ decomposes into $H$ and $H'$, and $e(N(s^i))\nearrow\infty$;
\item if $s^i=i/2$, each $N(s^i)$ decomposes into $H, H',$ and a fixed Seifert block,
and $e(N(s^i))\nearrow\infty$;
\end{itemize}
\item
let $N$ decompose into $k$ blocks, each homeomorphic to $P\times S^1$, with
one boundary component on $\partial N$, and the others glued cyclically 
(with some maps that do not match fibers). 
Take $s^i_j = i$ for all $i,j$. Then (for big $i\in\matN$) $N(s^i)$ is a \Sol\ torus fibering with 
$e(N(s^i))\nearrow\infty$.
\end{itemize}

\subsection{Beginning of the proof of Theorem \ref{convergence:teo}}
The rest of this Section is devoted to the proof of Theorem~\ref{convergence:teo}.
We set $s^i =(s_1^i,\ldots,s_k^i)$, and
we can suppose after passing to a subsequence 
that for each fixed $j$ the slopes $s_j^i$ are all distinct.

The following fact will be used below:
if $N(s^i)$ has a Seifert fibration with Euler number $e_i \nearrow\infty$,
we have $e(N(s^i))\nearrow\infty$, thanks to Proposition~\ref{small:Seifert:prop}.

\subsubsection{If $N$ is a solid torus} 
The foliation $\lambda$ is not the meridian $m$ by assumption. 
Proposition~\ref{distance:prop} implies that
$\Delta(s^i, m)\nearrow\infty$.
Therefore $N(s^i)$ is a lens space
with $e$ = $|\pi_1|=\Delta(s^i,m)\nearrow\infty$.

\subsubsection{If $N$ is a product $T\times[0,1]$} 
Both slopes $s^i_1$ and $s^i_2$ can be seen
inside $T$. Their limits $\lambda_1$ and $\lambda_2$ are distinct, 
otherwise they would bound an annulus.
Proposition~\ref{distance:prop}
implies that $\Delta(s_1^i, s_2^i)\nearrow\infty$,
hence $N(s^i)$ is again a lens space with $|\pi_1|\nearrow\infty$.

\subsubsection{If $N$ is a twisted bundle over the Klein bottle}
It fibers over the M\"obius strip. Its double cover is $A\times S^1$, where
$A$ is the annulus. The filled $N(s^i)$ is double
covered by $A\times S^1(s^i,-s^i)$, which is a lens space with $|\pi_1|=2|p^iq^i|$,
where $s^i = q^i/p^i$. Therefore $|\pi_1(N(q^i/p^i))| = |p^iq^i|\nearrow \infty$.

\subsubsection{If $N$ is another Seifert manifold} 
It fibers over an orbifold $\Sigma$ with $\chi(\Sigma)< 0$.
Also $N(s^i)$ fibers over an orbifold $\Sigma_i$, 
provided $s^i_j\neq\infty$ for all $j$, which
is certainly true on a subsequence. We have two cases:

\begin{itemize}
\item
if $\lambda_j\neq\infty$ for all $j$, we have 
$s_j^i=q_j^i/p_j^i$ with $|p_j^i|\nearrow\infty$ for all $j$.
Therefore $\chi(\Sigma_i)\searrow\chi(\Sigma)$ and hence
$\Vol(N(s^i))\nearrow\Vol(N)$;
\item
if $\lambda_j=\infty$ for some $j$, we have $\lambda_{j'}\neq\infty$ for all $j'\neq j$,
otherwise two foliations would coincide with the fiber of $N$, and bound an
annulus, contradicting the hypothesis. 
Therefore we get $|e(N(s^i))| = |e(N)+\sum_j s_j^i|\nearrow\infty$. 
\end{itemize}

\subsubsection{If $N$ is hyperbolic}
Thurston's Dehn
filling Theorem~\cite{bibbia} guarantees that $N(s^i)$ is hyperbolic except for
finitely many $i$'s, and that $\Vol(N(s^i))\nearrow\Vol(N)$.

\subsection{Induction on the number of geometric blocks}
We are left with the case where $N$ has a non-trivial 
geometric decomposition into blocks
$N_1,\ldots,N_t$. We proceed by induction on
$t$, having already considered the case $t=1$. 
Each $N_l$ is Seifert with $\chi<0$ or hyperbolic, 
by Proposition~\ref{geometric:prop}.

There is no compressing disc in $N$, and an incompressible annulus
is fibred in a single Seifert $N_l$.
Therefore our assumption
on $s$ may be replaced by the following:
\begin{equation} \label{assumption:eqn}
{\rm on\ each\ Seifert }\ N_l,\ {\rm at\ most\ one\ adjacent\ foliation\ }
\lambda_j {\rm\ is\ the\ fiber\ }\infty.
\end{equation}

\subsubsection{Persisting geometries}
Let $N_l(i)\subset N(s^i)$ be the connected submanifold consisting
of $N_l$ and all the filling solid tori adjacent to $N_l$.
Therefore $N(s^i) = \cup_{l=1}^t N_l(i)$. 
Using (\ref{assumption:eqn}) and Thurston's Dehn filling Theorem,
we can assume (after passing to a subsequence) 
that for each fixed $l$ one of the following holds:
\begin{itemize}
\item
$N_l$ is Seifert. No adjacent slope $s^i_j$ is the fiber $\infty$, and
at most one is an integer, for each $i$.
Moreover, one of the following holds:
\begin{itemize}
\item
$N_l(i)$ is Seifert with $\chi >0$, \emph{i.e.}~a solid torus, for all $i$;
\item
$N_l(i)$ is Seifert with $\chi = 0$, \emph{i.e.}~an interval bundle, for all $i$;
\item
$N_l(i)$ is Seifert with $\chi <0$, with the same fiber, for all $i$;
\end{itemize}
\item
$N_l$ and $N_l(i)$ are hyperbolic for all $i$, and $\partial N_l(i)$
has the same preferred slopes for all $i$.
\end{itemize}

\subsubsection{Solid tori}
We first consider the case
$N_l(i)$ is a solid torus for some $l$. 
Therefore $N_l$ fibers over an orbifold $\Sigma$ with $\chi(\Sigma)<0$, 
and $N_l(i)$ fibers over an orbifold $\Sigma_i$ with $\chi(\Sigma_i) >0$,
obtained by capping $\Sigma$.
Since $N_l$ has at most one integer filling,
it is easy to see that $\Sigma$ must be 
a pair-of-pants $P$ or an annulus with one cone point. In both cases,
$N_l(i)$ is obtained by filling $P \times S^1$ 
with two slopes $q_1^i\in\matZ$ and $q_2^i/p_2^i\in\matQ$ 
(where $q_2^i/p_2^i = q/p$ is fixed if $\Sigma$ is an annulus with cone point
$2\pi/p$).
The meridian of $N_l(i)$, read in the third boundary component of $P\times S^1$, 
is easily seen to be $-(q_1^i + q_2^i/p_2^i)$.

We have $|q_1^i|\nearrow\infty$ and $q_2^i/p_2^i\to q_2/p_2\neq\infty$, 
hence $|q_2^i/p_2^i|<K$.
Therefore the meridian $-(q_1^i + q_2^i/p_2^i)$
converges to the fiber $\infty$ of $N_l$.

Set $N_* = N\setminus N_l$. We modify $s^i$ to a set of slopes $s_*^i$
for $\partial N_*$, by removing the two slopes $q_1^i$ and $q_2^i/p_2^i$,
previously adjacent to $N_l$, and by adding
the meridian of $N_l(i)$. We get $N_*(s^i_*) = N(s^i)$.

By what is said above, the new meridian converges to
a slope $\lambda_*$, which is 
the fiber of the removed $N_l$. If the block $N_{l'}$ adjacent 
to $N_l$ is Seifert, $\lambda_*$ is not the fiber of $N_{l'}$ 
(because the fibers of $N_l$ and $N_{l'}$ do not match!).
Therefore assumption (\ref{assumption:eqn})
certainly holds also for $N_*$. Since $N_*$ has $t-1$ blocks, 
our main assertion holds for $N_*(s^i_*)$ 
by the induction hypothesis, and we are done.

\subsubsection{Products}
We are left with the case where no $N_l(i)$ is a solid torus.
It follows in particular that $N(s^i)$ is irreducible for all $i$
(because every $N_l(i)$ is $\partial$-irreducible). 

Some $N_l(i)$ are interval bundles.
We now prove that (for big $i\in\matN$) the blocks of 
the geometric decomposition of $N(s^i)$ are
the complements of these bundles.
By Proposition~\ref{geometric:prop}, it suffices to show
that if a sequence of product bundles connects a Seifert block $B$
with another Seifert block $B'$ (or a Klein bottle $K$) 
of the decomposition, the fiber of $B$ is not isotopic 
through the products to a fiber of $B'$ (or $K$).

By assumptions (1) and (2), one product $N_l(i)$ arise only
when $N_l = P\times S^1$ ($P$ is the pair-of-pants), 
and $N_l(i)$ is the filling of $N_l$ with some slope $q^i\in\matZ$. 
The same $N(s^i)$ is realized by removing the product $N_l(i)$
and gluing the adjacent blocks directly via a map $F_i$ that twists $q^i$ times
along the fiber of $N$. Since $|q^i|\to\infty$, we have
$\Delta(F_i(\eta),\eta')\nearrow\infty$ for any slopes $\eta, \eta'$ 
on the adjacent blocks distinct from the fiber of $N_l$.

Twisted bundles have a similar behaviour. If $N_l(i)$ is a twisted bundle,
$N_l$ is either an annulus with one point of cone angle $\pi$, filled with a 
half-integer $q^i/2$, or a 
M\"obius strip with a hole filled with an integer $q^i$.
In both cases we have $|q^i|\nearrow\infty$, 
whose effect is to twist the gluing
map with the adjacent block along the fiber $\gamma$.

It follows that the distance between the fibers of $B$ and $B'$ (or $K$)
goes to $\infty$. Hence the blocks of the geometric decomposition
of $N(s^i)$ consist of the complements of the interval bundles, and
the distance between some of their preferred slopes goes to $\infty$.
Therefore, if there is at least one interval bundle we get $e(N(s^i))\nearrow\infty$
and we are done.

Actually, if all the $N_l(i)$'s are interval bundles, then each $N(s^i)$ is a
\Sol-manifold, and the discussion above easily shows that
$e(N(s^i))\nearrow\infty$, as required.

\subsubsection{Geometric blocks}
We are now left with the case that there is no interval bundle.
Each $N_l(i)$ is then a geometric block of $N(s^i)$ by Proposition~\ref{geometric:prop}.

Suppose first some foliation $\lambda_j$ is the fiber $\infty$ of some Seifert block $N_l$.
We have $s^i_j = q^i/p^i \to \infty$, and $|s^i_{j'}|<K$ for every
other slope adjacent to $N_l$, by assumption (\ref{assumption:eqn}).
Therefore the sum $S^i$ of the $s^i_j$ adjacent to $N_l$ goes to $\infty$.
Let $F_i$ be a section of $N_l(i)$ (minus its singular fibers) 
for each $i$.
Since $S^i\to\infty$, on a subsequence
the component of $\partial F_i$ in a fixed boundary torus $T$
tends to the fiber $\gamma$ of $N_l(i)$.
The other block adjacent to $T$ has a preferred slope $\gamma'$ distinct from $\gamma$ and not
depending on $i$: 
then $\Delta(\gamma',\partial F_i)\nearrow\infty$ implies $e(N(s^i))\nearrow\infty$.

Finally, it remains to consider the case no foliation is 
a fiber of a Seifert block.
Let $N_l(i)$ be one block. If it is hyperbolic
we have $\Vol(N_l(i))\nearrow\Vol(N_l)$.
If it is Seifert, every slope $s^i_j=q^i/p^i$ converges to
a limit $\lambda_j\neq\infty$, hence $|p^i|\to\infty$ and 
we also have $\Vol(N_l(i))\nearrow\Vol(N_l)$.
Therefore we have $\Vol(N(s^i))\nearrow\Vol(N)$ and we are done.
The proof of Theorem~\ref{convergence:teo} is now complete.
\finedimo

\section{Links} \label{links:section}
We prove here Theorem~\ref{link:complements:teo} and Corollary~\ref{link:cor}.
We start by showing the following result.
A \emph{twist} along
an essential disc or annulus $S$ in a $3$-manifold $M$ is
a self-homeomorphism of $M$ constructed by cutting $M$ along $S$ and gluing it back after
a full twist.
\begin{teo}\label{self:teo}
Let $M$ be a $3$-manifold with boundary consisting of $k$ tori, and $N$ be a closed
$3$-manifold. There are only finitely many vectors of slopes $s=(s_1,\ldots, s_k)$ with
$M(s)=N$, that are not related by combinations of self-homeomorphisms of partial fillings
of $M$. 
The self-homeomorphisms consist of twists along essential 
discs or annuli.
\end{teo}

\begin{figure}[t]
\begin{minipage}{.02\textwidth}\hfil
\end{minipage}
\begin{minipage}{.2\textwidth}
\centering
\mettifig{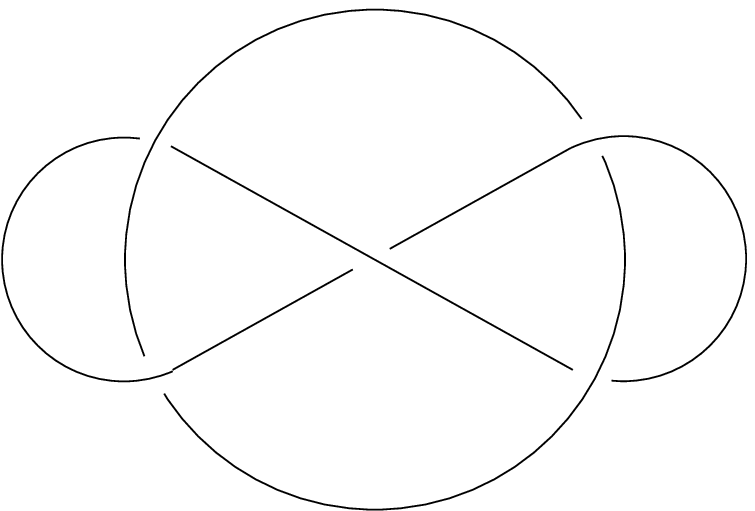, width = .9\textwidth}
\vglue -4mm
\end{minipage}
\begin{minipage}{.02\textwidth}\hfil
\end{minipage}
\begin{minipage}{.66\textwidth}
{The pairs of slopes yielding $S^3$ are $(\infty, 1/n)$ and $(1/n,\infty)$,
for $n\in\matZ$ \cite{magic}.
An $\infty$-filling on one component
gives a solid torus, and twisting along its meridian we transform every 
such pair into $(\infty, \infty)$.
}
\end{minipage}
\caption{The infinitely many fillings of the Whitehead link complement yielding the $3$-sphere
are all related by twists along discs in partial fillings. The partial fillings here
are solid tori. 
We use the standard meridian-longitude basis for link complements.}
\label{Whitehead:fig}
\end{figure}

Note that the discs and annuli are contained in the partial fillings of $M$, and not
necessarily in $M$, see Fig.~\ref{Whitehead:fig} as an example.

\subsubsection{Proof of Theorem~\ref{self:teo}}
It suffices to prove the theorem when $M$ is irreducible, because $M$ and $N$ 
split into summands in finitely many ways.
Suppose by contradiction we have infinitely many $s^i=(s_1^i,\ldots,s_k^i)$ with
$M(s^i)=N$, that are not related via combinations of twists along discs and annuli
in partial fillings. After passing to a subsequence, we can suppose either $s_j^i$ does
not depend on $i$, or $s_j^i\to\lambda_j$ with $s_j^i\neq\lambda_j$
for all $i$.
If $s_j^i$ does not depend on $i$ we can permanently fill $M$ with $s_j^i$,
and get a new $M'$ with less boundary components. Therefore we are left with the case that
$s^i\to\lambda=(\lambda_1,\ldots,\lambda_k)$ with $\lambda_j\neq s_j^i$ for all $i$.

Since $N=M(s^i)$, Theorem~\ref{convergence:teo} implies that some limits bound
a disc or an annulus. If the first
case holds, $M$ is a solid torus, and our assertion is easily proved.
Now suppose $\lambda_1$ and $\lambda_2$ bound an annulus. Set $\lambda_1 = \lambda_2 = \infty$.
By twisting along the annulus
we transform $(s^i_1, s^i_2)$ into $(s^i_1 \pm 1, s^i_2 \mp 1)$.
We can therefore suppose $|s^i_1|\leqslant 1$ for all $i$. 
If we get only finitely many slopes in $\{s^i_1\}_{i\in\matN}$ we proceed by induction on $k$ as above.
Otherwise, on a subsequence we have $s^i_1\to\lambda'_1\in [-1,1]$ 
distinct from the original limit $\lambda_1=\infty$.

We therefore get a new sequence of slopes, converging on one boundary torus to 
a different limit. Theorem~\ref{convergence:teo} applies
again, and we find another annulus. By iterating
this argument, we find at least two annuli
incident to distinct slopes on the same boundary component. 
Then $M$ is the product $T\times [0,1]$.

Finally, the theorem for $M=T\times [0,1]$ is easily proved: by twisting along annuli
we recover any self-homeomorphism of $M$ that fixes the boundary. Via such homeomorphisms
we transform each $s_1^i$ to a fixed slope $s$, and we fill permanently
$M$ along $s$, getting a solid torus, as above.
\finedimo

\subsubsection{Proof of Theorem~\ref{link:complements:teo}}
It follows from Theorem~\ref{self:teo}, thanks to the following 
1-1 bijective correspondence:
$$\left\{{\rm links\ in\ } M {\rm\ sharing\ the\ same\ complement\ } N
\right\}_{/_{{\rm Aut}(M)}} \ \longleftrightarrow \ \{ s\ |\ N(s) = M\}_{/_{{\rm Aut}(N)}}$$
which translates a twist modifying a link into a twist modifying some slopes.
\finedimo

\subsubsection{Proof of Corollary~\ref{link:cor}}
If $M=S^3$, homeomorphisms of $M$ need not to be considered because there are only
two up to isotopy. 
If there is an essential disc in $S^3\setminus \interior{N(K)}$ for some 
component $K$, then $K$ is the unknot and the disc is spanned by $K$.

Suppose now there is an essential annulus $A$ in the complement 
of $\interior{N(K_1\cup K_2)}$. If one component of $\partial A$ is
the meridian of $K_1$ or $K_2$, then $A$ extends to a disc in $S^3$ as above.
Otherwise, we can see $A$ spanning $K_1$ and $K_2$ inside $S^3$ (\emph{i.e.}~with
$\partial A = K_1\cup K_2)$), and
we must prove that $A$ winds only once along either $K_1$ or $K_2$. If not,
its regular neighborhood $N(A)$ in $S^3$ fibers over a disc with two cone points,
hence it is not a solid torus. 

Every torus in $S^3$ bounds a solid torus on one side.
Hence $S^3\setminus\interior{N(A)}$ is a solid torus, which extends the fibration of $N(A)$ to
a fibration of $S^3$ onto $S^2$. Since
$S^3$ does not fiber over $S^2$ with $3$ singular fibers, the new fiber is non-singular,
hence $S^3\setminus \interior{N(K_1\cup K_2)}$ is homeomorphic to $T\times [0,1]$.
Therefore $K_1\cup K_2$ is the Hopf link, and
a twist along $A$ is generated by a composition of twists along 
the two discs spanned by $K_1$ and $K_2$. 
\finedimo

\section{Two-handles and four-manifolds} \label{handle:slides:section}

We prove here Theorem~\ref{main:teo} and all the results about $4$-manifolds
stated in the Introduction.

\subsection{Dehn surgery}
We recall some well-known facts.
Let $L\subset M$ be a link in some closed oriented 3-manifold $M$. 
A \emph{Dehn surgery} on $L$ is a Dehn filling on the complement
$M_L = M\setminus\interior {N(L)}$. For each component $K$ of $L$, a \emph{longitude}
is an essential closed curve $l$ in $\partial N(K)$ 
which forms, together with the meridian $m$, a basis $(m,l)$ for 
$H_1(\partial N(K),\matZ)$. The choice of a longitude on each component $K$
of $L$ is a \emph{framing}, and it allows us to describe a slope $\pm(pm+ql)$ on some $K$ via
the number $p/q$.
When $L$ is framed, we can denote by $M_L(s)$ the manifold obtained by surgering
$L$ according to some vector $s=(s_1,\ldots,s_k)\in\matQ^k$.

\begin{teo} \label{links:teo}
Let $L\subset M$ be a framed link in 
some closed $3$-manifold $M$, with the following property:
there is no $2$-sphere $\Sigma\subset M$ intersecting only one or two components of $L$,
and each component in a single point.

For every closed $3$-manifold $N$ there is a constant $K>0$ 
(depending on $L\subset M$, its framing, and $N$) 
such that if $N=M_L(s_1,\ldots,s_k)$ 
then $|s_j|<K$ for some $j$.
\end{teo}
\begin{proof}
If $M_L$ is reducible,
then $L$ and $M$ split into summands, and we are left to prove 
the same theorem for each summand (because $N$ can split
into two summands in finitely many ways: as usual, if $M_L$ contains a non-separating
sphere we get a $S^2 \times S^1$-summand).
We can then restrict to the case $M_L$ is irreducible. 
We suppose by contradiction that there are infinitely many vectors
$s^i=(s^i_1,\ldots,s^i_k)\in\matQ^k$, $i\in\matN$,
with $M_L(s^i)=N$ for all $i$ 
and $s^i_j\to\infty$ for each $j$.
Therefore $s^i$ converges to the meridians 
$s = (\infty,\ldots,\infty)$ of $L$.

Since there are no spheres intersecting $L$ in one or
two points in distinct components,
there is no disc or annulus in $M_L$ bounded by the meridians.
Therefore Theorem \ref{convergence:teo} applies, and
a subsequence of $M_L(s^i)$ 
consists of distinct manifolds, distinguished
by their volume or Euler number: a contradiction.
\end{proof}

Note that a 2-sphere $\Sigma$ intersecting $L$ as stated must be non-separating.
Hence this condition is only needed when $M$ has some $S^2\times S^1$ summand.

\begin{rem}
Theorem~\ref{links:teo} is false without the hypothesis on $2$-spheres.
For instance, if $M=S^2\times S^1$ and $L=\{$pt$\}\times S^1$ we get $M_L(s)=S^3$
for every $s\in\matZ$.
\end{rem}

\subsection{Handle slides}
We refer to~\cite{GS} for the definition of handles and their main properties in the
$4$-manifolds setting.
In the statement of Theorem~\ref{main:teo}, let $L\subset M$ have $k$ components. 
After fixing an arbitrary framing on $L$, every other framing is
encoded by a vector of integers $s=(s_1,\ldots,s_k)$ colouring the components of $L$.
The result of attaching $2$-handles along that framing is a cobordism between $M$ and $N=M_L(s)$.

We can think about the $2$-handles as being attached to $M$ simultaneously, or one 
at each time, following some ordering of the components of $L$. In the latter case,
the $i$-th handle is attached along a knot $K_i$ contained in some manifold $M^i$,
giving rise to a cobordism between $M^i$ and $M^{i+1}$,
with $M^1=M$ and $M^{k+1} = N$. By isotoping the last knot $K_k$ inside $M^k$ we get the same
cobordisms, but we may change the initial framed link $L$. Such a modification of $L$ is
called a \emph{handle slide}\footnote{This definition is slightly more general than the usual one,
where one handle slides over another one.}.

We now use Theorem~\ref{links:teo} to prove the following stronger
version of Theorem~\ref{main:teo}.

\begin{teo} \label{main:refined:teo}
Let $M$ and $N$ be two closed $3$-manifolds, and $L\subset M$ be a link.
Among cobordisms between $M$ and $N$ obtained
by adding $2$-handles along $L$, only finitely many are
pairwise not related by combinations of handle slides.
\end{teo}
\begin{proof}
We do an induction on the number $k$ of components of $L$.
We fix a framing on $L$.
Suppose by contradiction we have infinitely many
cobordisms that are pairwise not related
by handle slides. 
Each cobordism is determined by a vector
of integers $s^i=(s_1^i,\ldots,s_k^i)\in\matZ^k$, $i\in\matN$.
We have $N=M_L(s^i)$ for all $i$.

If $L$ is a knot, Theorem~\ref{links:teo} implies 
that $L$ intersects a 2-sphere
$\Sigma$ in a point. We can visualize a neighborhood of $\Sigma$ in $M$ via a trivial
$0$-framed unknot encircling a portion of $L$
as in Fig.~\ref{slide1:fig}-left. 
The handle slide shown in Fig.~\ref{slide1:fig}
transforms $L$ into itself and changes the colour $a$ into $a+2$. 
A combination of slides (or of their inverses) transforms the colour
on $L$ into $0$ or $1$. Therefore there are at most two cobordisms
that are not related via handle slides. 

\begin{figure}
\begin{center}
\mettifig{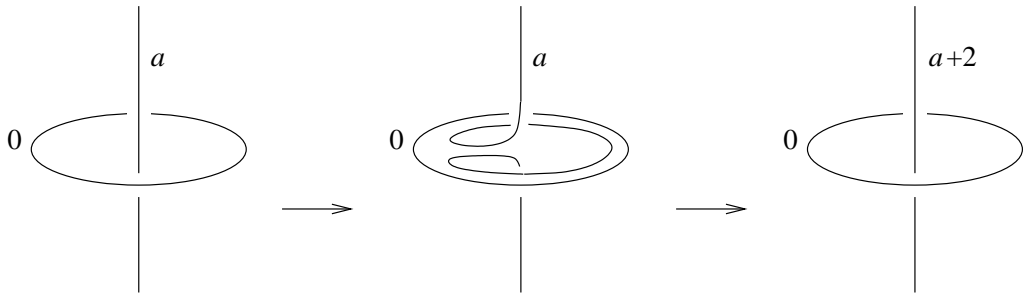} 
\caption{A handle slide, visualized as a Kirby move~\cite{K, GS} over a $0$-framed unknot.
We use the blackboard framing.}
\label{slide1:fig}
\end{center}
\end{figure}

Let now $L$ have some $k\geqslant 2$ components, and
let $\calS_j=\{s_{j}^i\}_{i\in\matN}$ be the set
of integers colouring the $j$-th component of $L$.
We first consider the case $\calS_j$ contains infinitely many numbers for all $j$. 
Therefore we have
$\sup_{i\in\matN}|s_j^i| = \infty$ for all $j$. Theorem~\ref{links:teo}
then implies that there is a 2-sphere $\Sigma$ intersecting $L$ 
in one or two components,
and in a single point on each component. As above, we visualize $\Sigma$
via a $0$-framed unknot encircling a portion of $L$
as in Fig.~\ref{slide1:fig}-left or Fig.~\ref{slide2:fig}-top-left.
Let the $j$-th component of $L$
be one intersecting $\Sigma$.
Performing sufficiently many times the handle slides shown in
Figg.~\ref{slide1:fig} or~\ref{slide2:fig} (or their inverses) we can transform
$s^i_j$ into $0$ or $1$ for 
each $i\in\matN$. Therefore $\calS_{j}$ becomes finite.

\begin{figure}
\begin{center}
\mettifig{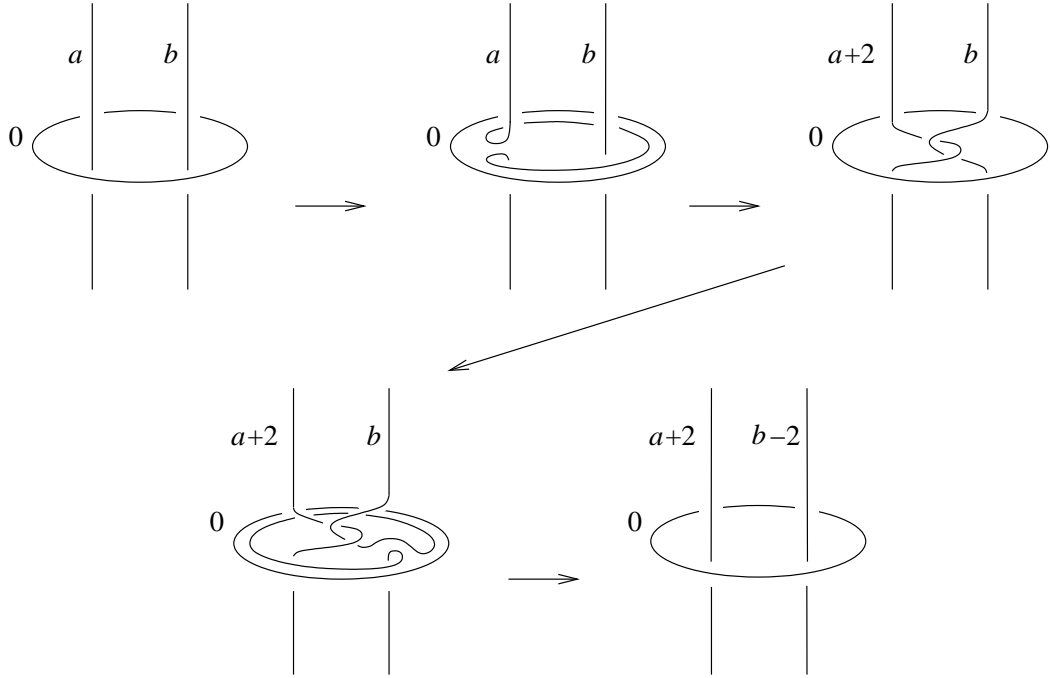} 
\caption{A combination of handle slides.}
\label{slide2:fig}
\end{center}
\end{figure}

We are left to consider the case $\calS_j$ is finite for some $j$. 
We restrict ourselves to a subsequence  
of $i\in\matN$ so that $s_j^i$ is the same integer for all $i$.
We can sort the handles so that the $j$-th becomes the first one.
This handle gives a cobordism between $M=M^1$ and a fixed $M^2$, not depending on $i$.
The other $k-1$ handles are attached to some fixed link $L^2\subset M^2$ and
give cobordisms between $M^2$ and $N$ that depend on $i$.
By our induction hypothesis some of these cobordisms
are related by a combination of handle slides. The handle
slides on $L^2\subset M^2$ translate to 
handle slides on $L\subset M$, and we are done. 
\end{proof}

\subsection{Kirby diagrams}
We now prove all the assertions stated in Subsection~\ref{closed:subsection}.
Combining Theorem~\ref{main:teo} with 
Laudenbach-Poenaru Theorem~\cite{LP}, we get the following.

\begin{teo}\label{handles:teo}
The set of closed orientable 
$4$-manifolds constructed with one zero-handle, $h$ one-handles,
some $2$-handles attached along a fixed 
link $L\subset \#_hS^2\times S^1$, and some $3$- and $4$-handles, is finite.
\end{teo}
\begin{proof}
The $2$-handles form a cobordism between 
$M=\#_hS^2\times S^1$ and some $N=\#_{h'}S^2\times S^1$,
with $h' =$ rk$(H_1(N)) \leqslant$ rk$(H_1(M\setminus L))$ bounded above.
By Theorem~\ref{main:teo}, we have a finite number of 
such cobordisms for each $h'$.
By Laudenbach-Poenaru Theorem~\cite{LP}, each such cobordism can be filled with $0$-, $1$-handles 
and $3$-, $4$-handles in a unique way.
\end{proof}

\subsubsection{Proof of Corollary~\ref{Kirby:cor}}
Only finitely many link diagrams have weight at most $n$.
A diagram with $h$ couples of $2$-spheres describes a fixed link 
in $\#_hS^2\times S^1$, and yields 
only finitely many manifolds by Theorem~\ref{handles:teo}.
\finedimo

\subsection{Locally flat $2$-skeleta} \label{shadows:subsection}
We prove here Theorem~\ref{polyhedra:teo}.
Let $P\subset M$ be a locally flat $2$-skeleton.
We first show that local flatness can be extended to the $1$-skeleton $P^1$.
That is, a regular neighborhood $N(P^1)$ is contained in a smooth $3$-submanifold of $M$.

\begin{figure}
\begin{center}
\mettifig{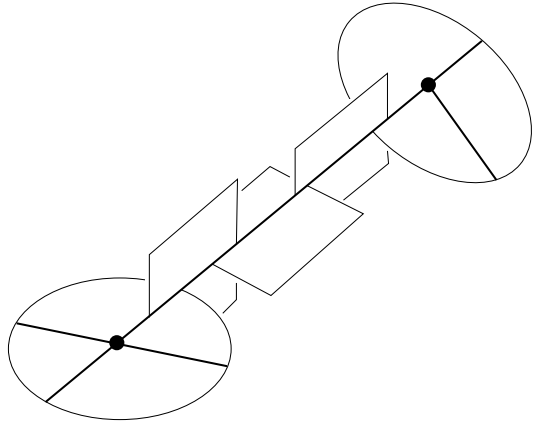, width = 7cm} 
\caption{Extending local flatness to the one-skeleton: here $P^2$ is a graph in $\matR^3$,
locally contained in $2$-discs that have to be patched together.}
\label{poly:fig}
\end{center}
\end{figure}

Take a small $3$-ball containing each vertex of $P^1$. A neighborhood of each edge $e$ 
is covered by a finite number of open smooth $3$-balls
(see Fig.~\ref{poly:fig}-left for a smaller-dimensional picture)
and we want to patch them together.
Take a sufficiently thin tubular neighborhood of $e$, which we identify
with $B^3\times [0,1]$, where $e$ becomes $\{0\}\times [0,1]$.
We can suppose
every $3$-ball is the zero-set of a smooth function $f_i: B^3\times I_i \to\matR$,
where $\{I_i\}$ are open intervals covering $[0,1]$.
We have $\nabla f_i\neq 0$ everywhere on $B^3\times I_i$, and
$P\cap (B^3\times I_i)$ is contained in the zero-set of $f_i$.

We now patch the
functions $f_i$ and $f_{i+1}$ on a thinner tubular 
neighborhood of $\{0\}\times (I_i\cup I_{i+1})$.
Fix a point $y_0\in I_i\cap I_{i+1}$. 
Take a smooth function 
$g:B^3\times [0,1]\to [0,1]$ with $g(x,y)= 0$ for $y<y_0-\epsilon$ and
$g(x,y)=1$ for $y>y_0+\epsilon,$ where
$\epsilon$ is a sufficiently small number.
We patch $f_i$ and $f_{i+1}$ to the function:
$$h = (1-g)f_i + gf_{i+1},$$ defined over $B^3\times (I_i\cup I_{i+1})$.
The polyhedron $P\cap \big(B^3\times (I_i\cup I_{i+1})\big)$ is
contained in the zero-set of $h$. 

It remains to prove that $\nabla h\neq 0$ on
$\{0\}\times (I_i\cup I_{i+1})$: this implies that $\nabla h\neq 0$ 
on a thin tubular neighborhood, and we are done.
On $\{0\}\times (I_i\cup I_{i+1})$ we have:
$$\nabla h = (1-g)\nabla f_i + g\nabla f_{i+1}$$
because $f_i$ and $f_{i+1}$ are always zero. Then $\nabla h$
equals $\nabla f_i$ or $\nabla f_{i+1}$, which are non-zero, everywhere 
except on $\{0\}\times [y_0-\epsilon, y_0+\epsilon]$. 

Up to changing sign to $f_i$,
we can suppose $\nabla f_i$ and $\nabla f_{i+1}$ are not opposite
on $y_0$, and thus on the whole $\{0\}\times [y_0-\epsilon, y_0+\epsilon]$
by taking $\epsilon$ small. This implies that
$\nabla h\neq 0$ also there. 

We have proved that a regular neighborhood $N(P^1)$ of 
the $1$-skeleton $P^1$ of $P\subset M^4$ is contained in
a $3$-dimensional smooth manifold. Therefore there is a
(possibly non-orientable) $3$-dimensional handlebody $H^3\subset M^4$ intersecting
$P$ in $N(P^1)$, which thickens in $M^4$ to a $4$-dimensional handlebody $H^4$,
uniquely determined as the orientable line bundle over $H^3$.
Let $L\subset\partial H^3\subset\partial H^4$ be the link given by $L=P\cap\partial H^3$.
The $H^4$ is made of $0$- and $1$-handles.
The faces $P\setminus N^1(P)$ of $P$ thicken in $M$ to $2$-handles, attached along $L$.
The complement is made of $3$- and $4$-handles.

We now prove that only finitely many pairs $(H^3, L)$, and thus $(H^4,L)$,
can be obtained from a fixed $P$ (and varying $M$). Together
with Theorem~\ref{handles:teo}, this implies that only finitely many $M$'s can be obtained
from a fixed $P$. 
Actually, there are only finitely many pairs $(H^3, N^1(P))$, since in general a
$2$-dimensional polyhedron has
only finitely many $3$-dimensional thickenings.
Each such thickening is constructed as follows:
embed each vertex $v$ in $D^3$ as a cone on a graph on $\partial D^3$
(there are only finitely many choices, corresponding to embeddings of the link of $v$ inside
$S^2$). The thickening of $N(P^1)$ can be then extended along each edge in at most
two ways.
Therefore there are only finitely many possible $(H^3,L)$ at the end.
\finedimo

\begin{rem}\label{PL:rem}
A $2$-dimensional polyhedron that has locally flat faces is not necessarily
locally flat, \emph{i.e.}~locally flatness of the $1$-skeleton is a serious hypothesis.
For instance, take a closed braid $L$ in $S^2\times S^1$.
Let $P\subset D^3\times S^1$ be the polyhedron defined by making cones over the braid,
\emph{i.e.}~every slice $P\cap \big(D^3\times\{{\rm pt}\}\big)$ is the cone over the points
$L\cap \big(S^2\times\{{\rm pt}\}\big)$ with center $\{0\}\times\{{\rm pt}\}$.
The proof of Theorem~\ref{handles:teo} shows that if $P$ is
locally flat then it is contained in a properly embedded solid torus, and hence
$L$ is contained in a torus $T\subset S^2\times S^1$.
Therefore if $L$ is sufficiently complicated
the polyhedron $P$ is not locally flat, not even after an isotopy.
\end{rem}
\begin{rem}
It is also essential that every point has a neighborhood in a \emph{smooth} $3$-ball.
For instance, for generic braids $L$
every point in the polyhedron $P$ constructed in Remark~\ref{PL:rem}
has a neighborhood
contained in a PL $3$-ball (\emph{i.e.}~it is locally PL-flat!), constructed
as follows: if the point lies inside a face we are done. Otherwise
it is some $\{0\}\times \{x_0\}$, and by genericity of $L$ there is at most one triple
of points in $S^2\times \{x_0\}$ lying inside a circle of maximum length,
and no such triple in $S^2\times \{x\}$ for all points $x$ near $x_0$.
It is easy to deduce that there is a PL polygon $\gamma_x\subset S^2\times \{x\}$
having the points $L\cap (S^2\times \{x\})$ as vertices for all $x$ close to $x_0$, 
which has a PL-dependence on $x$.
The cone of $\gamma_x$ in $D^3\times\{x\}$ is
a PL $2$-disc. The union of such cones is a PL $3$-disc.
(Note that these PL $3$-discs can be patched together precisely when the polygon
$\gamma_x$ runs continuously over all $x\in S^1$, 
\emph{i.e.}~when $L$ is contained in a torus!)
\end{rem}

\subsection{Rescalings}
\begin{figure}
\begin{center}
\mettifig{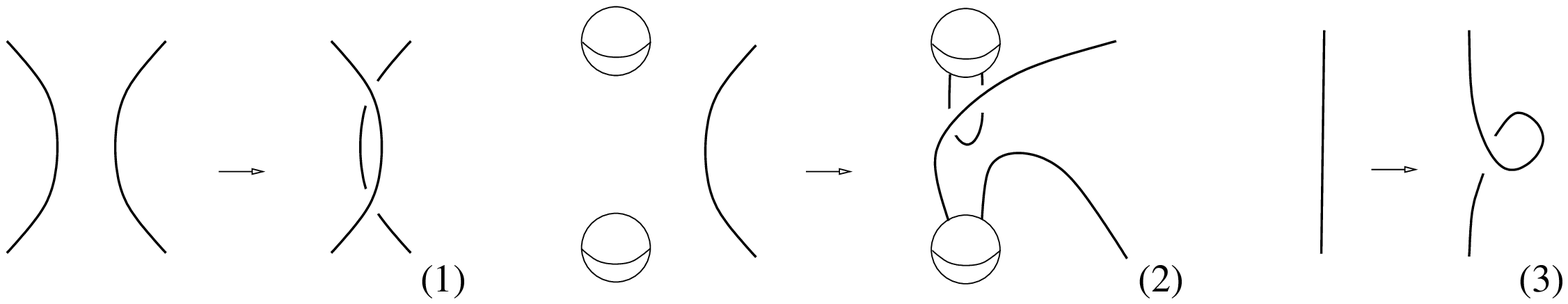, width = 14 cm} 
\caption{Sliding strands we can connect the whole diagram.}
\label{moves:fig}
\end{center}
\end{figure}

We now prove Theorem~\ref{rescalings:teo}. 
We first show $U_n\subset V_{3n}$. If $M\in U_n$, there is
a Kirby diagram for $M$ of weight at most $n$.
If the diagram is not connected, we can connect components 
via Reidemeister moves as in Fig.~\ref{moves:fig}-(1), producing $2$ new
crossings for each move. If the diagram contains a couple of discs not attached to any strand,
we slide a strand over it as in Fig.~\ref{moves:fig}-(2). (If there is no
strand at all we add a $1$-framed trivial knot.) 
Finally, if the diagram has no crossings we produce one via a Reidemeister move as in 
Fig.~\ref{moves:fig}-(3).

\begin{figure}
\begin{center}
\mettifig{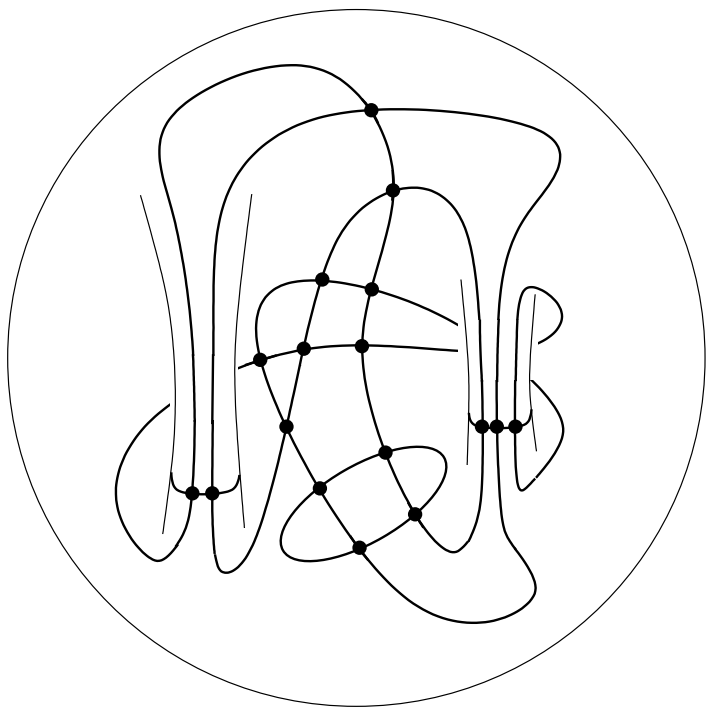} 
\caption{Passing from the Kirby diagram of Fig.~\ref{Kirby:fig} 
to a shadow: first construct a surface with genus $k$ and
$k$ cocore discs $\Sigma$ (here $k=2$), then add the core discs of the $2$-handles to it, following
the diagram. Dots represent vertices of the shadow.}
\label{diagram2shadow_new:fig}
\end{center}
\end{figure}

The resulting diagram is connected, has some crossings, and has weight at most $3n$.
It contains some $k$ couples of discs, denoting $1$-handles.
We now construct a shadow representing the same $4$-manifold.
Take first a big disc containing the whole diagram, and close it to a $2$-sphere.
Then, for each couple of discs in the diagram, add a $1$-handle 
and a cocore disc as in Fig.~\ref{diagram2shadow_new:fig}.
The resulting polyhedron $\Sigma$ is a genus-$k$ surface with $k$ discs attached, and it is a
shadow of the submanifold consisting of all $0$- and $1$-handles (\emph{i.e.}~the 
submanifold is a neighborhood of $\Sigma$ plus a $3$-handle).

Now we attach a disc for each $2$-handle following the diagram, and we get our shadow $P$.
Each crossing and each non-closed strand produces a vertex for $P$: therefore
$P$ has $3n$ vertices at most. Moreover, $P$ is easily seen to be special,
because the original diagram is connected and contains some crossings.

\begin{figure}
\begin{center}
\mettifig{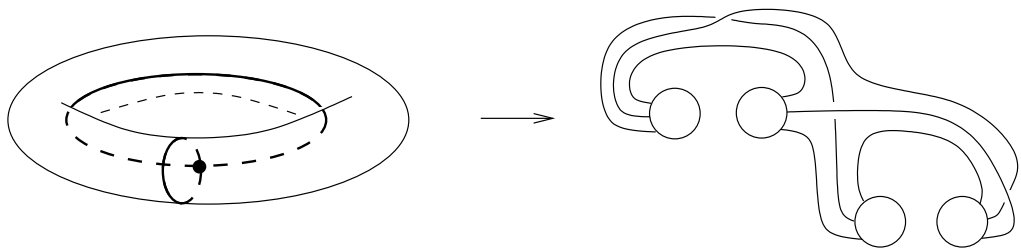} 
\caption{Passing from a shadow to a Kirby diagram. The integers on the strands
will depend on the embedding of the shadow in the $4$-manifold
(\emph{i.e.}~from the \emph{gleams} of the shadow~\cite{Tu}).}
\label{diagram2shadow:fig}
\end{center}
\end{figure}

Concerning $V_n\subset U_{9n+8}$, there is a well-known translation of a shadow $P$ into
a Kirby diagram~\cite{CoTh}: 
first, take a diagram representing $P$ on the plane~\cite{Mat}. Then
cut the diagram in $n+1$ points, corresponding to $n+1$ edges of $P$ 
having a tree as a complement, and add a pair of discs at each cut,
as in Fig.~\ref{diagram2shadow:fig}.
We get one crossing for each vertex and
at most $3$ crossings for each of the $n+1$ cut edges, hence $4n+3$ crossings at most.
There are precisely $3(n+1)$ strands and $2(n+1)$ discs. Therefore the weight of the Kirby diagram is
at most $9n+8$.
\finedimo

\subsection{Fundamental groups}
We prove here Theorem~\ref{groups:teo}. We start with a similar
result concerning special polyhedra.

\begin{prop} \label{groups:prop}
The number of groups that are fundamental groups of
special polyhedra with $n$ vertices grows as $n^n$.
\end{prop}
\begin{proof}
It is shown in~\cite{FriMaPe3} that there are at most $n^{C\cdot n}$ special
polyhedra with $n$ vertices, and at least $n^{c\cdot n}$ of them
are spines of distinct hyperbolic $3$-manifolds with geodesic boundary,
for some $0<c<C$ and all $n\gg 0$.
By Mostow rigidity, distinct hyperbolic 
manifolds have distinct fundamental groups (see~\cite{Fri} for a proof of this
fact in the geodesic boundary case). 
Therefore the $n^{c\cdot n}$ polyhedra
have distinct fundamental groups.
\end{proof}

\subsubsection{Proof of Theorem~\ref{groups:teo}}
Proposition~\ref{groups:prop} implies that the number of fundamental groups
of compact $4$-manifolds with boundary having a special shadow with $n$ vertices
grows as $n^n$. To recover the same result for closed $4$-manifolds, it suffices
to use Theorem~\ref{rescalings:teo} to pass from shadows to diagrams and vice versa,
together with the following fact: any diagram of weight $n$ represents a $4$-manifold $M$
made of $0$-, $1$-, and $2$-handles.
By encircling every strand with a small $0$-framed unknot we get a diagram for
its double $DM$~\cite{GS} of weight at most $4n$. The double is closed, and 
$\pi_1(DM)=\pi_1(M)$.
\finedimo

\subsection{Simply connected manifolds}
We now prove Propositions~\ref{simply:connected:prop} and~\ref{homeomorphic:prop}.

\subsubsection{Proof of Proposition~\ref{simply:connected:prop}}
Take $k$ unknots, coloured with $\pm 1$ and representing some
$\#_h\matCP^2\#_{k-h}\overline{\matCP^2}$. 
With these $k$ unknots we get $k/2$ distinct manifolds at least. 
With $k$ ranging from $1$ to $n$ we get $n^2/4$ manifolds at least.

\subsubsection{Proof of Proposition~\ref{homeomorphic:prop}}
By Freedman's Theorem~\cite{F} simply connected closed
smooth $4$-manifolds are determined up to homeomorphism by their intersection forms.
We now compute how many manifolds can have an intersection form of fixed rank $n$.
If the form is odd, it is of type $k\langle 1 \rangle \oplus h\langle -1\rangle$ with $k+h=n$,
and we suppose $k-h\geqslant 0$ up to switching the orientation 
(the indefinite case is a general result on unimodular forms~\cite{GS}, while the definite case $h=0$
follows from Donaldson~\cite{D}). Then we get $n/2+1$ forms at most.
If the form is even, it is of type $2kE_8\oplus lH$ by Rohlin's Theorem~\cite{Roh, GS}, and we get
$n/16+1$ forms at most. Summing up, there are at most $\frac 9{16}n+2$ manifolds up to
homeomorphism.

A diagram with weight $n$ produces a bilinear form of rank at most $n$. Therefore
$U_n$ contains at most $\sum_{i=1}^n (\frac 9{16}i+2) = \frac 9{32}n^2 + o(n^2) < \frac{5}{16} n^2$
manifolds up to homeomorphism, when $n\gg 0$.
\finedimo

\end{document}